\documentclass[12pt]{amsart}
\usepackage{amssymb,latexsym}

\def\AA{\mathcal{A}}
\def\cc{\mathbf{c}}
\def\FF{\mathcal{F}}
\def\QQ{\mathbb{Q}}
\def\SS{\mathcal{S}}
\def\xx{\mathbf{x}}
\def\XX{\mathbf{X}}
\def\ex{\mathbf{ex}}
\def\ZZ{\mathbb{Z}}

\newtheorem{theorem}{Theorem}[section]

\theoremstyle{definition}

\newtheorem{definition}[theorem]{Definition}
\newtheorem{example}[theorem]{Example}

\begin{document}

\title{Quantum cluster algebras: Oberwolfach talk, February 2005}

\author{Andrei Zelevinsky}
\address{\noindent Department of Mathematics, Northeastern University,
  Boston, MA 02115, USA}
\email{andrei@neu.edu}

\begin{abstract}
This is an extended abstract of my talk at the
Oberwolfach-Workshop ``Representation Theory of Finite-Dimensional
Algebras" (February 6 - 12, 2005).
It gives self-contained and simplified definitions of quantum
cluster algebras introduced and studied in a joint work with
A.~Berenstein (math.QA/0404446).
\end{abstract}

\date{February 13, 2005}

\thanks{Research supported in part
by NSF (DMS) grant \# 0200299.}

\maketitle

\makeatletter
\renewcommand{\@evenhead}{\tiny \thepage \hfill  A.~BERENSTEIN and  A.~ZELEVINSKY \hfill}

\renewcommand{\@oddhead}{\tiny \hfill QUANTUM CLUSTER ALGEBRAS
 \hfill \thepage}
\makeatother

Cluster algebras were introduced and studied by S.~Fomin and
A.~Zelevinsky in \cite{ca1,ca2,ca3}. This is a family of
commutative rings designed to serve as an algebraic framework for
the theory of total positivity and canonical bases in semisimple
groups and their quantum analogs. Here we report on a joint work
with A.~Berenstein \cite{qca}, where we introduce and study
quantum deformations of cluster algebras.

We start by recalling the definition of cluster algebras (of
geometric type). Let $m$ and $n$ be two positive integers with $m
\geq n$. Let $\FF$ be the field of rational functions over $\QQ$
in $m$ independent (commuting) variables.

\begin{definition}
\label{def:seed} A \emph{seed} in $\FF$ is a pair $(\tilde \xx,
\tilde B)$, where
\begin{itemize}

\item $\tilde \xx = \{x_1, \dots, x_m\}$
is a free (i.e., algebraically independent) generating set
for~$\FF$.

\item $\tilde B$ is an $m \times n$ integer matrix with rows
labeled by $[1,m] = \{1, \dots, m\}$ and columns labeled by an
$n$-element subset $\ex \subset [1,m]$, such that, for some
positive integers $d_j \,\, (j \in \ex)$, we have $d_i b_{ij} = -
d_j b_{ji}$ for all $i, j \in \ex$.
\end{itemize}
The subset $\xx = \{x_j: j \in \ex\} \subset \tilde \xx$ (resp.
$\cc = \tilde \xx - \xx$) is called the \emph{cluster} (resp. the
\emph{coefficient set}) of a seed $(\tilde \xx, \tilde B)$. The
seeds are defined up to a relabeling of elements of $\tilde \xx$
together with the corresponding relabeling of rows and columns of
$\tilde B$.
\end{definition}

\begin{definition}
\label{def:seed-mutation} Let $(\tilde \xx, \tilde B)$ be a seed
in $\FF$. For any $k \in \ex$, the \emph{seed mutation} in
direction~$k$ transforms $(\tilde \xx, \tilde B)$ into a seed
$(\tilde \xx', \tilde B')$ given by:
\begin{itemize}
\item
$\tilde \xx' = \tilde \xx - \{x_k\} \cup \{x'_k\}$, where $x'_k
\in \FF$ is determined by the \emph{exchange relation}
\begin{equation}
\label{eq:exchange-rel-xx} x'_k = x_k^{-1} \,(\prod_{\substack{i
\in [1,m] \\ b_{ik}>0}} x_i^{b_{ik}} + \prod_{\substack{i \in
[1,m] \\ b_{ik}<0}} x_i^{-b_{ik}}) \ .
\end{equation}

\item
The entries of~$\tilde B'$ are given by
\begin{equation}
\label{eq:matrix-mutation} b'_{ij} =
\begin{cases}
-b_{ij} & \text{if $i=k$ or $j=k$;} \\[.05in] b_{ij} +
\displaystyle\frac{|b_{ik}| b_{kj} + b_{ik} |b_{kj}|}{2} &
\text{otherwise.}
\end{cases}
\end{equation}

\end{itemize}
\end{definition}

The seed mutations generate an equivalence relation: we say that
two seeds $(\tilde \xx, \tilde B)$ and $(\tilde \xx', \tilde B')$
are \emph{mutation-equivalent} if $(\tilde \xx', \tilde B')$ can
be obtained from $(\tilde \xx, \tilde B)$ by a sequence of seed
mutations.

Fix a mutation-equivalence class~$\SS$ of seeds. Let $\mathcal{X}
\subset \FF$ denote the union of clusters, and~$\cc$ the common
coefficient set of all seeds from~$\SS$. The \emph{cluster
algebra}~$\AA(\SS)$ associated with~$\SS$ is the $\ZZ[\cc^{\pm
1}]$-subalgebra of~$\FF$ generated by~$\mathcal{X}$.

We now define a family of $q$-deformations of~$\AA(\SS)$. The
following setup is a simplified version of that in \cite{qca}. The
main idea is to deform each extended cluster $\tilde \xx$ to a
quasi-commuting family $\tilde \XX = \{X_1, \dots, X_m\}$
satisfying
\begin{equation}
\label{eq:Xi-q-com} X_i X_j = q^{\lambda_{ij}} X_j X_i
\end{equation}
for some skew-symmetric integer $m \times m$ matrix $\Lambda =
(\lambda_{ij})$. Let $\FF_q$ denote the skew-field of fractions of
the ring $\ZZ[q^{\pm 1/2}, X_1, \dots, X_m]$, where $X_1, \dots,
X_m$ are algebraically independent variables satisfying
\eqref{eq:Xi-q-com}. For any $a = (a_1, \dots, a_m) \in \ZZ^m$, we
set
\begin{equation}
\label{eq:M-Xi} X^a = q^{\frac{1}{2} \sum_{i > j}\lambda_{ij} a_i
a_j} X_1^{a_1} \cdots X_m^{a_m} \ .
\end{equation}

\begin{definition}
\label{def:q-free-gen-set} A \emph{free generating set} for
$\FF_q$ is a subset $\{Y_1, \dots, Y_m\} \subset \FF_q$ of the
following form: $Y_j = \varphi(X^{c_j})$, where $\varphi$ is a
$\QQ(q^{\pm 1/2})$-linear automorphism of $\FF_q$, and $\{c_1,
\dots, c_m\}$ is a basis of the lattice~$\ZZ^m$.
\end{definition}

Note that the subset $\{Y_1, \dots, Y_m\}$ can be used instead of
$\{X_1, \dots, X_m\}$ in the definition of the ambient
field~$\FF_q$, with the matrix~$\Lambda$ replaced by $C^T \Lambda
C$, where~$C$ is the matrix with columns $c_1, \dots, c_m$.

\begin{definition}
\label{def:q-seed} A \emph{quantum seed} in $\FF_q$ is a pair
$(\tilde \XX, \tilde B)$, where
\begin{itemize}

\item $\tilde \XX = \{X_1, \dots, X_m\}$
is a free generating set for~$\FF_q$.

\item $\tilde B$ is a $m \times n$ integer matrix with rows
labeled by $[1,m]$ and columns labeled by an $n$-element subset
$\ex \subset [1,m]$, which is \emph{compatible} with the
matrix~$\Lambda$ given by \eqref{eq:Xi-q-com}, in the following
sense: for some positive integers $d_j \,\, (j \in \ex)$, we have
\begin{equation}
\label{eq:compatibility} \sum_{k = 1}^m b_{kj} \lambda_{ki} =
\delta_{ij} d_j \quad (j \in \ex, \, i \in [1,m]) \ .
\end{equation}
\end{itemize}
As in Definition~\ref{def:seed}, the quantum seeds are defined up
to a relabeling of elements of $\tilde \XX$ together with the
corresponding relabeling of rows and columns of $\tilde B$.
\end{definition}

Note that \eqref{eq:compatibility} implies that $d_i b_{ij} = -
d_j b_{ji}$ for all $i, j \in \ex$, i.e., $\tilde B$ is as in
Definition~\ref{def:seed}.

\begin{example}
Let $m = 2n$, $\ex = [1,n]$, and let $\tilde B$ be of the form
$$\tilde B = \left(\!\!\begin{array}{c} B \\ I \\
\end{array}\!\!\right) \, ,$$
where $I$ is the identity $n \times n$ matrix. Here $B$ is an
arbitrary integer $n \times n$ matrix satisfying $d_i b_{ij} = -
d_j b_{ji}$ for some positive integers $d_1, \dots, d_n$: in other
words, $B$ is \emph{skew-symmetrizable}, that is, $DB$ is
skew-symmetric, where $D$ is the diagonal matrix with diagonal
entries  $d_1, \dots, d_n$. An easy calculation shows that the
skew-symmetric matrices~$\Lambda$ compatible with~$\tilde B$ in
the sense of \eqref{eq:compatibility} are those of the form
\begin{equation}
\label{eq:lambda-principal} \Lambda = \left(\!\!\begin{array}{cc}
\Lambda_0 & -D - \Lambda_0 B \\ D - B^T \Lambda_0 & -DB + B^T
\Lambda_0 B\\
\end{array}\!\!\right) \, ,
\end{equation}
where $\Lambda_0$ is an arbitrary skew-symmetric integer $n \times
n$ matrix.
\end{example}

\begin{definition}
\label{def:q-seed-mutation} Let $(\tilde \XX, \tilde B)$ be a
quantum seed in $\FF_q$. For any $k \in \ex$, the \emph{quantum
seed mutation} in direction~$k$ transforms $(\tilde \XX, \tilde
B)$ into a quantum seed $(\tilde \XX', \tilde B')$ given by:
\begin{itemize}
\item
$\tilde \XX' = \tilde \XX - \{X_k\} \cup \{X'_k\}$, where $X'_k
\in \FF_q$ is given by
\begin{equation}
\label{eq:q-exchange-rel-xx} X'_k = X^{- e_k + \sum_{b_{ik} > 0}
b_{ik} e_i} \ + \
 X^{- e_k - \sum_{b_{ik} < 0} b_{ik} e_i} \ ,
\end{equation}
where the terms on the right are defined via \eqref{eq:M-Xi}, and
$\{e_1, \dots, e_m\}$ is the standard basis in $\ZZ^m$.
\item
The matrix entries of $\tilde B'$ are given by
\eqref{eq:matrix-mutation}.
\end{itemize}
\end{definition}

The fact that $(\tilde \XX', \tilde B')$ is a quantum seed is not
automatic: for the proof see \cite[Proposition~4.7]{qca}.

Based on definitions~\eqref{def:q-seed} and
\eqref{def:q-seed-mutation}, one defines the \emph{quantum cluster
algebra} associated with a mutation-equivalence class of quantum
seeds, in exactly the same way as the ordinary cluster algebra. It
is shown in \cite{qca} that practically all the structural results
on cluster algebras obtained in \cite{ca1,ca2,ca3} extend to the
quantum setting. This includes the Laurent phenomenon obtained in
\cite{ca1,fz-laurent,ca3} and the classification of cluster
algebras of finite type given in \cite{ca2}.

\end{document}